# A Simple Local Variational Iteration Method and Related Algorithm for Nonlinear Science and Engineering

Xuechuan Wang[*], Qiuyi Xu, Satya N. Atluri

Department of Mechanical Engineering
Center for Advance Research in Engineering Sciences
Texas Tech University, Lubbock, TX, 79415, USA

## Abstract

A very simple and efficient local variational iteration method for solving problems of nonlinear science is proposed in this paper. The analytical iteration formula of this method is derived first using a general form of first order nonlinear differential equations, followed by straightforward discretization using Chebyshev polynomials and collocation method. The resulting numerical algorithm is very concise and easy to use, only involving highly sparse matrix operations of addition and multiplication, and no inversion of the Jacobian in nonlinear problems. Apart from the simple yet efficient iteration formula, another extraordinary feature of LVIM is that in each local domain, all the collocation nodes participate in the calculation simultaneously, thus each local domain can be regarded as one "node" in calculation through GPU acceleration and parallel processing. For illustration, the proposed algorithm of LVIM is applied to various nonlinear problems including Blasius equations in fluid mechanics, buckled bar equations in solid mechanics, the Chandrasekhar equation in astrophysics, the low-Earth-orbit equation in orbital mechanics, etc. Using the built-in highly optimized *ode45* function of MATLAB as a comparison, it is found that the LVIM is not only very accurate, but also much faster by an order of magnitude than *ode45* in all the numerical examples, especially when the nonlinear terms are very complicated and difficult to evaluate.

## 1. Introduction

A general numerical method for solving nonlinear systems is fundamental to investigate and predict nonlinear phenomena of various models in nonlinear science. So far, the most widely used and simplest method for solving nonlinear systems, often modelled by partial differential or ordinary differential equations, is the direct finite difference method [1]. Based on Taylor series approximation, it is often applied as central difference, forward difference, and backward difference schemes. Another broad class of methods are based on weighted residual approximation [2], including finite element methods [3] and finite volume methods [4]. They are often applied to solve complicated problems in solid and fluid mechanics. In solving nonlinear ordinary differential equations, the finite difference method is more often preferred for its simplicity and possibility to easily obtain high order approximation. Here we list six nonlinear equations in fluid mechanics, structural mechanics, and astrophysics, etc.

---

[*] Corresponding author. Email address: xc.wang@ttu.edu

Blasius equation:

$$\frac{\partial u}{\partial x}+\frac{\partial v}{\partial y}=0, \quad u\frac{\partial u}{\partial x}+v\frac{\partial u}{\partial y}=\nu\frac{\partial^2 u}{\partial y^2}.$$

Emden-Chandrasekhar equation:

$$\frac{1}{\xi^2}\frac{d}{d\xi}\left(\xi^2\frac{d\psi}{d\xi}\right)=e^{-\psi}.$$

Mathieu equation:

$$\frac{d^2 x}{dt^2}+(\delta-\varepsilon\cos t)x=0.$$

Pendulum equation:

$$\frac{d^2\theta}{dt^2}+\frac{g}{l}\sin\theta=0.$$

Buckled bar equations:

$$EI\frac{d^2\theta}{ds^2}=-P\sin\theta \quad \text{(vertical dead load)},$$

$$EI\frac{d^2\theta}{ds^2}=-P\cos(\theta-\alpha)\sin\theta \quad \text{(perpendicular follower load)},$$

$$EI\frac{d^2\theta}{ds^2}=-P\sin(\theta-\alpha)\sin\theta \quad \text{(tangent follower load)}.$$

Euler's elastica equation:

$$\frac{dy}{dx}=\frac{a^2-c^2+x^2}{\sqrt{(c^2-x^2)(2a^2-c^2+x^2)}}.$$

These equations can of course be solved using finite difference methods such as the well-known Runge-Kutta (RK) method, Hilber-Hughes-Taylor (HHT)-$\alpha$ method [5], and Newmark-$\beta$ method [6]. The best performance of these methods is superseded by the *ode* functions built in MATLAB, particularly the *ode45* function in solving low dimensional systems. However, although *ode* functions are sufficient to solve low-dimensional systems, we aim to propose a superior method and its related algorithm that can achiever better performance in computational efficiency and accuracy, not just in simple problems, but more significantly in problems with high dimensions and complicated calculations.

Asymptotic methods have been regularly used to analyze nonlinear systems with "weak" nonlinearity [7]. Among them, the variational iteration method (VIM) [8], the Adomian decomposition method (ADM) [9] and the Picard iteration method (PIM) [10] have received wide attention. However, these methods involve too many symbolic calculations even for simple problems. In practice, it is often impossible to utilize them to obtain accurate solutions on computers. Considering that, we propose a Local Variational Iteration Method (LVIM), where the whole domain is divided into small segments and the solution is obtained in each local domain. Unlike the prementioned asymptotic methods, the generalized Lagrange multiplier [11] in LVIM is optimally determined so that it can achieve best performance in the correcting process. Moreover, the initial guess of LVIM can take a very simple form. Unlike the conventional VIM, there is no need to pay much attention to the construction of an initial guess. A local linear function could be good enough, as long as it meets the boundary condition of the problem. The LVIM can be used to predict long-term transient motion, even in complex nonlinear problems with chaos.

The finite difference methods mentioned above have already been well popularized in both academic research and commercial software, but computational difficulties still exist for large scale systems with complex nonlinear terms, especially when high accuracy and computational efficiency are required. In this paper, we introduce a method based on variational principle and collocation discretization in local domain, named Local Variational Iteration Method (LVIM). The applications of LVIM in nonlinear science include but not limit to $n$-body problem in celestial mechanics [12], molecular dynamics [13], high-dimensional nonlinear structural dynamics [14], and other problems in fluid mechanics, solid mechanics, heat transfers, etc. The $n$-body problem in celestial mechanics can be formulated as

$$m_i \frac{d^2 \mathbf{q}_i}{dt^2} = \sum_{\substack{j=1 \\ j \neq i}}^{n} \frac{G m_i m_j (\mathbf{q}_j - \mathbf{q}_i)}{\|\mathbf{q}_j - \mathbf{q}_i\|^3} = -\frac{\partial U}{\partial \mathbf{q}_i},$$

where $i = 1, 2, 3, ..., n$, $m_i$ are point mass, $\mathbf{q}_i$ are position vectors, and $U$ is the total potential energy. Although the formulation is quite simple, it is very tricky to solve when dimensions and complexity of gravitational force grows. It is known that for $n \geq 3$, the system become chaotic, meaning that the motion of celestial bodies is usually not as determined as the periodic motion in two-body systems. This brings extra computational difficulty because the system itself is not computationally stable, and a small computational error could be significantly enlarged in a relatively short time. Moreover, if high-fidelity models are necessary, the perturbation terms in the $n$-body system will be very time-consuming to calculate. Take a low Earth orbit for example, there could be a thousand to a million perturbation terms in addition to the central gravity force. The commonly used finite difference methods such as Gauss-Jackson method, high-order RK method are clumsy in such cases [15]. To solve such problems, one needs a method with low computational complexity, fast convergence rate, and eligibility for parallel computation. Another similar example is the molecular dynamical problem [16]. It can be formulated as

$$m_i \frac{d^2 \mathbf{r}_i}{dt^2} = \mathbf{F}_i = -\frac{\partial U}{\partial \mathbf{r}_i} \ , \ U(\mathbf{r}_{ij}) = 4\varepsilon \left[ \left(\frac{\sigma}{\mathbf{r}_{ij}}\right)^{12} - \left(\frac{\sigma}{\mathbf{r}_{ij}}\right)^{6} \right],$$

where $U(\mathbf{r}_1, \mathbf{r}_2, ..., \mathbf{r}_N)$ is the potential function describing the potential energy of a system of $N$ atoms. The use of potential energy to describe the quantum interaction among atoms is justified through Born-Oppenheimer approximation [17] that assumes the electrons adjust to new atomic positions much faster than the motion of the atomic nuclei. The potential function we presented herein is the most commonly seen Lennard-Jones model [18], which is one of pair potentials. In structural and material dynamics, the governing equations are often discretized into ordinary differential equations in the following form

$$\mathbf{M}\frac{d^2\mathbf{x}}{dt^2} + \mathbf{C}\frac{d\mathbf{x}}{dt} + \mathbf{N}(\frac{d\mathbf{x}}{dt}, \mathbf{x}) = \mathbf{F}(t),$$

where $\mathbf{M}$, $\mathbf{C}$, $\mathbf{N}$ are mass, damping, and nonlinear stiffness matrices, $\mathbf{F}$ is the external load. Even for some simple structures, this system could well exceed thousands of variables and dimensions. Using conventional methods in literature, the computation could be very time consuming and inaccurate. Moreover, for stiff systems, the numerical calculation of conventional methods can easily diverge due to ill conditioning [19].

In this paper, the variational iteration formula of LVIM is derived first, from which we found that two classical asymptotic methods for solving nonlinear systems, the Picard iteration and the Adomian decomposition method, can be regarded as its particular versions [20]. This is elucidated by introducing the concept of general Lagrange multiplier. Then, the collocation discretization is made to the variational iteration formula in local domain, so that the formula can be carried out numerically on computers. The resulting algorithm of LVIM is very simple and computationally very efficient. All the matrix operations in it are sparse and easy to achieve since it only involves addition and multiplication, but no inversion of Jacobian matrix in nonlinear problems. The formula of LVIM is highly modular, meaning that the user need not to derive the formula for different problems, but only need to replace the modules related with the particular problem, i.e. the governing equation and the Jacobian matrix. By using orthogonal basis function and appropriate discretized nodes, LVIM can achieve very high accuracy and efficiency in solving various problems in nonlinear science. For illustration, seven problems are selected as numerical examples, i.e. the Blasius equation in boundary layer theory [21], the Emden-Chandrasekhar [22] and Chandrasekhar's white dwarf equation [23], the Mathieu equation [24], the large-amplitude pendulum equation [25], the buckled bar equation [26], the Euler's elastica equation [27], and the low-Earth-orbit equation. The numerical results are compared with those of the highly optimized *ode45* built in MATLAB. Discussions are made in both numerical and physical aspects for each problem. We show that LVIM is at least an order of magnitude superior in accuracy and speed as compared to the highly optimized function *ode45* in MATLAB.

## 2. Methodology

The derivation of Local Variational Iteration Method (LVIM) involves two stages: the construction of functional iterative formula of LVIM and the numerical discretization for computational convenience. These two stages of derivation are presented as following. The ready-to-use formula of LVIM is given in Eq. (26).

## 2.1 Construction of functional iterative formula of LVIM

Generally, for solving a system of first order nonlinear differential equations

$$\frac{d\mathbf{x}}{d\tau} = \mathbf{g}(\mathbf{x},\tau), \ \tau \in [t_0,t], \ t \in [t_0,t_f], \tag{1}$$

the LVIM method approximates the solution at any time $t$ with an initial approximation $\mathbf{x}_0(\tau)$ and the correctional iterative formula as

$$\mathbf{x}_{n+1}(t) = \mathbf{x}_n(\tau)\big|_{\tau=t} + \int_{t_0}^{t} \boldsymbol{\lambda}(\tau)\{\dot{\mathbf{x}}_n(\tau) - \mathbf{g}[\mathbf{x}_n(\tau),\tau]\}d\tau, \tag{2}$$

where $\boldsymbol{\lambda}(\tau)$ is a matrix of Lagrange multipliers yet to be determined. Eq. (2) indicates that the $(n+1)$th correction to the analytical solution $\mathbf{x}_{n+1}$ involves the addition of $\mathbf{x}_n$ and a feedback weighted optimal error in the solution $\mathbf{x}_n$ up to the current time $t$. $\boldsymbol{\lambda}(\tau)$ can be optimally determined by making the right-hand side of Eq. (2) stationary about $\delta\mathbf{x}_n(\tau)$, thus resulting the constraints.

$$\begin{cases} \delta\mathbf{x}_n(\tau)\big|_{\tau=t} : \mathbf{I} + \boldsymbol{\lambda}(\tau)\big|_{\tau=t} = \mathbf{0} \\ \delta\mathbf{x}_n(\tau) : \dot{\boldsymbol{\lambda}}(\tau) + \boldsymbol{\lambda}(\tau)\mathbf{J}(\tau) = \mathbf{0} \end{cases}, \ t_0 \leq \tau \leq t, \tag{3}$$

where $\mathbf{I}$ is a unit matrix, and $\mathbf{J}(\tau) = \partial\mathbf{g}(\mathbf{x}_n,\tau)/\partial\mathbf{x}_n$.

First order Taylor series approximation of $\boldsymbol{\lambda}(\tau)$ can be readily obtained from Eq. (3), in the following form

$$\boldsymbol{\lambda}(\tau) \approx -\mathbf{I} + \mathbf{J}(t)(\tau - t), \tag{4}$$

By substituting Eq. (4) into Eq. (2), and for simplicity assuming $\mathbf{G}(\tau) = \dot{\mathbf{x}}_n(\tau) - \mathbf{g}[\mathbf{x}_n(\tau),\tau]$, we have

$$\mathbf{x}_{n+1}(t) = \mathbf{x}_n(t) - \{\mathbf{I} + \mathbf{J}(t)t\}\int_{t_0}^{t} \mathbf{G}(\tau)d\tau + \mathbf{J}(t)\int_{t_0}^{t} \tau\mathbf{G}(\tau)d\tau, \tag{7}$$

Two other ways of approximating $\boldsymbol{\lambda}(\tau)$ have been proposed by the authors leading to 3 distinct LVIM algorithms [28], but herein we just focus on using Eq. (7).

## 2.2 Numerical Discretization Through Collocation Method

The numerical discretization of Eq. (7) is achieved through collocation in time domain. Collocating $t_1, t_2, ..., t_M$ as time nodes, Eq. (7) leads to

$$\mathbf{x}_{n+1}(t_1) = \mathbf{x}_n(t_1) - \{\mathbf{I} + \mathbf{J}(t_1)t_1\}\int_{t_0}^{t_1}\mathbf{G}(\tau)d\tau + \mathbf{J}(t_1)\int_{t_0}^{t_1}\tau\mathbf{G}(\tau)d\tau$$

$$\mathbf{x}_{n+1}(t_2) = \mathbf{x}_n(t_2) - \{\mathbf{I} + \mathbf{J}(t_2)t_2\}\int_{t_0}^{t_2}\mathbf{G}(\tau)d\tau + \mathbf{J}(t_2)\int_{t_0}^{t_2}\tau\mathbf{G}(\tau)d\tau \quad (8)$$

$$\vdots$$

$$\mathbf{x}_{n+1}(t_M) = \mathbf{x}_n(t_M) - \{\mathbf{I} + \mathbf{J}(t_M)t_M\}\int_{t_0}^{t_M}\mathbf{G}(\tau)d\tau + \mathbf{J}(t_M)\int_{t_0}^{t_M}\tau\mathbf{G}(\tau)d\tau$$

It is rewritten in matrix form as

$$\begin{bmatrix} \mathbf{x}_{n+1}(t_1) \\ \mathbf{x}_{n+1}(t_2) \\ \vdots \\ \mathbf{x}_{n+1}(t_M) \end{bmatrix} = \begin{bmatrix} \mathbf{x}_n(t_1) \\ \mathbf{x}_n(t_2) \\ \vdots \\ \mathbf{x}_n(t_M) \end{bmatrix} - (\tilde{\mathbf{I}} + \tilde{\mathbf{J}}\tilde{\mathbf{T}})\begin{bmatrix} \int_{t_0}^{t_1}\mathbf{G}(\tau)d\tau \\ \int_{t_0}^{t_2}\mathbf{G}(\tau)d\tau \\ \vdots \\ \int_{t_0}^{t_M}\mathbf{G}(\tau)d\tau \end{bmatrix} + \tilde{\mathbf{J}}\begin{bmatrix} \int_{t_0}^{t_1}\tau\mathbf{G}(\tau)d\tau \\ \int_{t_0}^{t_2}\tau\mathbf{G}(\tau)d\tau \\ \vdots \\ \int_{t_0}^{t_M}\tau\mathbf{G}(\tau)d\tau \end{bmatrix}, \quad (9)$$

where $\tilde{\mathbf{I}}$, $\tilde{\mathbf{J}}$, $\tilde{\mathbf{T}}$ are defined as

$$\tilde{\mathbf{I}} = \begin{bmatrix} \mathbf{I} & & & \\ & \mathbf{I} & & \\ & & \ddots & \\ & & & \mathbf{I} \end{bmatrix}, \tilde{\mathbf{J}} = \begin{bmatrix} \mathbf{J}(t_1) & & & \\ & \mathbf{J}(t_2) & & \\ & & \ddots & \\ & & & \mathbf{J}(t_M) \end{bmatrix}, \tilde{\mathbf{T}} = \begin{bmatrix} t_1\mathbf{I} & & & \\ & t_2\mathbf{I} & & \\ & & \ddots & \\ & & & t_M\mathbf{I} \end{bmatrix}, \quad (10)$$

Generally, $\mathbf{x}(t) = [x_1(t), x_2(t), ..., x_d(t), ... x_D(t)]^T$ is a D-dimensional vector of time-varying variables. We suppose $x_d(t)$, $d = 1, 2, ..., D$, are approximated by the linear combinations of $N$ orthogonal basis functions $\phi_{d,nb}(t)$,

$$x_d(t) = \sum_{nb=1}^{N} \alpha_{d,nb}\phi_{d,nb}(t) = \mathbf{\Phi}_d(t)\mathbf{A}_d, \quad (11)$$

where $\mathbf{\Phi}_d(t) = [\phi_{d,1}(t), \phi_{d,2}(t), ..., \phi_{d,N}(t)]$, $\mathbf{A}_d = [\alpha_{d,1}, \alpha_{d,2}, ..., \alpha_{d,N}]^T$. Note that $\mathbf{A}_d$ are unknown constant coefficients of basis functions. In practice, it is much easier to use the same set of orthogonal basis functions throughout the coding, thus we will simply denote $\mathbf{\Phi}_d(t)$ as $\mathbf{\Phi}(t)$. Obviously, at time nodes $t_1, t_2, ..., t_M$, Eq. (11) leads to

$$\begin{bmatrix} x_d(t_1) \\ x_d(t_2) \\ \vdots \\ x_d(t_M) \end{bmatrix} = \begin{bmatrix} \mathbf{\Phi}(t_1) \\ \mathbf{\Phi}(t_2) \\ \vdots \\ \mathbf{\Phi}(t_M) \end{bmatrix} \mathbf{A}_d, \text{ and } \begin{bmatrix} \dot{x}_d(t_1) \\ \dot{x}_d(t_2) \\ \vdots \\ \dot{x}_d(t_M) \end{bmatrix} = \begin{bmatrix} \dot{\mathbf{\Phi}}(t_1) \\ \dot{\mathbf{\Phi}}(t_2) \\ \vdots \\ \dot{\mathbf{\Phi}}(t_M) \end{bmatrix} \mathbf{A}_d. \quad (12)$$

It can be derived that

$$\begin{bmatrix} \dot{x}_d(t_1) \\ \dot{x}_d(t_2) \\ \vdots \\ \dot{x}_d(t_M) \end{bmatrix} = \begin{bmatrix} \dot{\Phi}(t_1) \\ \dot{\Phi}(t_2) \\ \vdots \\ \dot{\Phi}(t_M) \end{bmatrix} \begin{bmatrix} \Phi(t_1) \\ \Phi(t_2) \\ \vdots \\ \Phi(t_M) \end{bmatrix}^{-1} \begin{bmatrix} x_d(t_1) \\ x_d(t_2) \\ \vdots \\ x_d(t_M) \end{bmatrix}. \quad (13)$$

For simplicity, we define $\mathbf{Q}$, $\boldsymbol{\chi}_d$, $\dot{\boldsymbol{\chi}}_d$, and row rearranging matrix $\tilde{\mathbf{R}}$ so that

$$\mathbf{Q} = \begin{bmatrix} \dot{\Phi}(t_1) \\ \dot{\Phi}(t_2) \\ \vdots \\ \dot{\Phi}(t_M) \end{bmatrix} \begin{bmatrix} \Phi(t_1) \\ \Phi(t_2) \\ \vdots \\ \Phi(t_M) \end{bmatrix}^{-1}, \quad \boldsymbol{\chi}_d = \begin{bmatrix} x_d(t_1) \\ x_d(t_2) \\ \vdots \\ x_d(t_M) \end{bmatrix}, \quad \dot{\boldsymbol{\chi}}_d = \begin{bmatrix} \dot{x}_d(t_1) \\ \dot{x}_d(t_2) \\ \vdots \\ \dot{x}_d(t_M) \end{bmatrix}, \text{ and } \begin{bmatrix} \mathbf{x}(t_1) \\ \mathbf{x}(t_2) \\ \vdots \\ \mathbf{x}(t_M) \end{bmatrix} = \tilde{\mathbf{R}} \begin{bmatrix} \boldsymbol{\chi}_1 \\ \boldsymbol{\chi}_2 \\ \vdots \\ \boldsymbol{\chi}_D \end{bmatrix}. \quad (14)$$

In this way, we have

$$\begin{bmatrix} \dot{\mathbf{x}}(t_1) \\ \dot{\mathbf{x}}(t_2) \\ \vdots \\ \dot{\mathbf{x}}(t_M) \end{bmatrix} = \tilde{\mathbf{R}} \begin{bmatrix} \dot{\boldsymbol{\chi}}_1 \\ \dot{\boldsymbol{\chi}}_2 \\ \vdots \\ \dot{\boldsymbol{\chi}}_D \end{bmatrix} = \tilde{\mathbf{R}} \begin{bmatrix} \mathbf{Q} & & & \\ & \mathbf{Q} & & \\ & & \ddots & \\ & & & \mathbf{Q} \end{bmatrix} \begin{bmatrix} \boldsymbol{\chi}_1 \\ \boldsymbol{\chi}_2 \\ \vdots \\ \boldsymbol{\chi}_D \end{bmatrix} = \tilde{\mathbf{R}}\tilde{\mathbf{Q}} \begin{bmatrix} \boldsymbol{\chi}_1 \\ \boldsymbol{\chi}_2 \\ \vdots \\ \boldsymbol{\chi}_D \end{bmatrix}. \quad (15)$$

Thus,

$$\begin{bmatrix} \mathbf{G}(t_1) \\ \mathbf{G}(t_2) \\ \vdots \\ \mathbf{G}(t_M) \end{bmatrix} = \begin{bmatrix} \dot{\mathbf{x}}_n(t_1) \\ \dot{\mathbf{x}}_n(t_2) \\ \vdots \\ \dot{\mathbf{x}}_n(t_M) \end{bmatrix} - \begin{bmatrix} \mathbf{g}(t_1) \\ \mathbf{g}(t_2) \\ \vdots \\ \mathbf{g}(t_M) \end{bmatrix} = \tilde{\mathbf{R}}\tilde{\mathbf{Q}}\tilde{\mathbf{R}}^{-1} \begin{bmatrix} \mathbf{x}_n(t_1) \\ \mathbf{x}_n(t_2) \\ \vdots \\ \mathbf{x}_n(t_M) \end{bmatrix} - \begin{bmatrix} \mathbf{g}(t_1) \\ \mathbf{g}(t_2) \\ \vdots \\ \mathbf{g}(t_M) \end{bmatrix}. \quad (16)$$

Generally, $\mathbf{G}(t) = [G_1(t), G_2(t), ..., G_d(t), ... G_D(t)]^T$ is a D-dimensional vector having the same size as $\mathbf{x}(t)$. We approximate $G_d(t)$, $d = 1, 2, ..., D$, through interpolation using a set of orthogonal basis functions,

$$G_d(t) = \sum_{nb=1}^{N} \beta_{d,nb} \phi_{nb}(t) = \mathbf{\Phi}(t)\mathbf{B}_d, \quad (17)$$

where $\mathbf{\Phi}(t) = [\phi_1(t), \phi_2(t), ..., \phi_N(t)]$, $\mathbf{B}_d = [\beta_{d,1}, \beta_{d,2}, ..., \beta_{d,N}]^T$. Note that the same symbol $\mathbf{\Phi}(t)$ is used to denote the basis functions of both $x_d(t)$ and $G_d(t)$, because they can share the same set of basis functions. At time nodes $t_1, t_2, ..., t_M$, Eq. (17) leads to

$$\begin{bmatrix} G_d(t_1) \\ G_d(t_2) \\ \vdots \\ G_d(t_M) \end{bmatrix} = \begin{bmatrix} \mathbf{\Phi}(t_1) \\ \mathbf{\Phi}(t_2) \\ \vdots \\ \mathbf{\Phi}(t_M) \end{bmatrix} \mathbf{B}_d, \text{ and } \begin{bmatrix} \int_{t_0}^{t_1} G_d(\tau)d\tau \\ \int_{t_0}^{t_2} G_d(\tau)d\tau \\ \vdots \\ \int_{t_0}^{t_M} G_d(\tau)d\tau \end{bmatrix} = \begin{bmatrix} \int_{t_0}^{t_1} \mathbf{\Phi}(\tau)d\tau \\ \int_{t_0}^{t_2} \mathbf{\Phi}(\tau)d\tau \\ \vdots \\ \int_{t_0}^{t_M} \mathbf{\Phi}(\tau)d\tau \end{bmatrix} \mathbf{B}_d. \tag{18}$$

It can be derived that

$$\begin{bmatrix} \int_{t_0}^{t_1} G_d(\tau)d\tau \\ \int_{t_0}^{t_2} G_d(\tau)d\tau \\ \vdots \\ \int_{t_0}^{t_M} G_d(\tau)d\tau \end{bmatrix} = \begin{bmatrix} \int_{t_0}^{t_1} \mathbf{\Phi}(\tau)d\tau \\ \int_{t_0}^{t_2} \mathbf{\Phi}(\tau)d\tau \\ \vdots \\ \int_{t_0}^{t_M} \mathbf{\Phi}(\tau)d\tau \end{bmatrix} \begin{bmatrix} \mathbf{\Phi}(t_1) \\ \mathbf{\Phi}(t_2) \\ \vdots \\ \mathbf{\Phi}(t_M) \end{bmatrix}^{-1} \begin{bmatrix} G_d(t_1) \\ G_d(t_2) \\ \vdots \\ G_d(t_M) \end{bmatrix} \tag{19}$$

For simplicity, we define $\mathbf{P}$, $\mathbf{\Gamma}_d$, and $\int \mathbf{\Gamma}_d$ as

$$\mathbf{P} = \begin{bmatrix} \int_{t_0}^{t_1} \mathbf{\Phi}(\tau)d\tau \\ \int_{t_0}^{t_2} \mathbf{\Phi}(\tau)d\tau \\ \vdots \\ \int_{t_0}^{t_M} \mathbf{\Phi}(\tau)d\tau \end{bmatrix} \begin{bmatrix} \mathbf{\Phi}(t_1) \\ \mathbf{\Phi}(t_2) \\ \vdots \\ \mathbf{\Phi}(t_M) \end{bmatrix}^{-1}, \mathbf{\Gamma}_d = \begin{bmatrix} G_d(t_1) \\ G_d(t_2) \\ \vdots \\ G_d(t_M) \end{bmatrix}, \text{ and } \int \mathbf{\Gamma}_d = \begin{bmatrix} \int_{t_0}^{t_1} G_d(\tau)d\tau \\ \int_{t_0}^{t_2} G_d(\tau)d\tau \\ \vdots \\ \int_{t_0}^{t_M} G_d(\tau)d\tau \end{bmatrix}. \tag{20}$$

Then the integration terms in Eq. (9) are obtained as

$$\begin{bmatrix} \int_{t_0}^{t_1} \mathbf{G}(\tau)d\tau \\ \int_{t_0}^{t_2} \mathbf{G}(\tau)d\tau \\ \vdots \\ \int_{t_0}^{t_M} \mathbf{G}(\tau)d\tau \end{bmatrix} = \tilde{\mathbf{R}} \begin{bmatrix} \int \mathbf{\Gamma}_1 \\ \int \mathbf{\Gamma}_2 \\ \vdots \\ \int \mathbf{\Gamma}_D \end{bmatrix} = \tilde{\mathbf{R}} \begin{bmatrix} \mathbf{P} & & & \\ & \mathbf{P} & & \\ & & \ddots & \\ & & & \mathbf{P} \end{bmatrix} \begin{bmatrix} \mathbf{\Gamma}_1 \\ \mathbf{\Gamma}_2 \\ \vdots \\ \mathbf{\Gamma}_D \end{bmatrix} = \tilde{\mathbf{R}}\tilde{\mathbf{P}}\tilde{\mathbf{R}}^{-1} \begin{bmatrix} \mathbf{G}(t_1) \\ \mathbf{G}(t_2) \\ \vdots \\ \mathbf{G}(t_M) \end{bmatrix} \tag{21}$$

Similarly, we have

$$\begin{bmatrix} \int_{t_0}^{t_1} \tau \mathbf{G}(\tau) d\tau \\ \int_{t_0}^{t_2} \tau \mathbf{G}(\tau) d\tau \\ \vdots \\ \int_{t_0}^{t_M} \tau \mathbf{G}(\tau) d\tau \end{bmatrix} = \tilde{\mathbf{R}} \tilde{\mathbf{P}} \tilde{\mathbf{R}}^{-1} \begin{bmatrix} t_1 \mathbf{G}(t_1) \\ t_2 \mathbf{G}(t_2) \\ \vdots \\ t_M \mathbf{G}(t_M) \end{bmatrix}. \tag{22}$$

Substituting Eq. (16) into Eqs. (21, 22), and then Eqs. (21, 22) into Eq. (9), it is obtained that

$$\begin{bmatrix} \mathbf{x}_{n+1}(t_1) \\ \mathbf{x}_{n+1}(t_2) \\ \vdots \\ \mathbf{x}_{n+1}(t_M) \end{bmatrix} = \begin{bmatrix} \mathbf{x}_n(t_1) \\ \mathbf{x}_n(t_2) \\ \vdots \\ \mathbf{x}_n(t_M) \end{bmatrix} + \left\{ \tilde{\mathbf{J}} \tilde{\mathbf{R}} \tilde{\mathbf{P}} \tilde{\mathbf{R}}^{-1} \tilde{\mathbf{T}} - \left( \tilde{\mathbf{I}} + \tilde{\mathbf{J}} \tilde{\mathbf{T}} \right) \tilde{\mathbf{R}} \tilde{\mathbf{P}} \tilde{\mathbf{R}}^{-1} \right\} \left\{ \tilde{\mathbf{R}} \tilde{\mathbf{Q}} \tilde{\mathbf{R}}^{-1} \begin{bmatrix} \mathbf{x}_n(t_1) \\ \mathbf{x}_n(t_2) \\ \vdots \\ \mathbf{x}_n(t_M) \end{bmatrix} - \begin{bmatrix} \mathbf{g}(t_1) \\ \mathbf{g}(t_2) \\ \vdots \\ \mathbf{g}(t_M) \end{bmatrix} \right\} \tag{23}$$

Since $\tilde{\mathbf{R}}$ is simply row rearranging matrix, we denote $\tilde{\mathbf{P}} = \tilde{\mathbf{R}} \tilde{\mathbf{P}} \tilde{\mathbf{R}}^{-1}$ and $\tilde{\mathbf{Q}} = \tilde{\mathbf{R}} \tilde{\mathbf{Q}} \tilde{\mathbf{R}}^{-1}$. Eq. (23) becomes

$$\begin{bmatrix} \mathbf{x}_{n+1}(t_1) \\ \mathbf{x}_{n+1}(t_2) \\ \vdots \\ \mathbf{x}_{n+1}(t_M) \end{bmatrix} = \begin{bmatrix} \mathbf{x}_n(t_1) \\ \mathbf{x}_n(t_2) \\ \vdots \\ \mathbf{x}_n(t_M) \end{bmatrix} + \left\{ \tilde{\mathbf{J}} \left( \tilde{\mathbf{P}} \tilde{\mathbf{T}} - \tilde{\mathbf{T}} \tilde{\mathbf{P}} \right) - \tilde{\mathbf{P}} \right\} \left\{ \tilde{\mathbf{Q}} \begin{bmatrix} \mathbf{x}_n(t_1) \\ \mathbf{x}_n(t_2) \\ \vdots \\ \mathbf{x}_n(t_M) \end{bmatrix} - \begin{bmatrix} \mathbf{g}(t_1) \\ \mathbf{g}(t_2) \\ \vdots \\ \mathbf{g}(t_M) \end{bmatrix} \right\} \tag{24}$$

For convenience of coding, we denote constant matrix $\tilde{\mathbf{P}} \tilde{\mathbf{T}} - \tilde{\mathbf{T}} \tilde{\mathbf{P}}$ as $\tilde{\mathbf{H}}$, and

$$\tilde{\mathbf{x}} = \begin{bmatrix} \mathbf{x}(t_1) \\ \mathbf{x}(t_2) \\ \vdots \\ \mathbf{x}(t_M) \end{bmatrix}, \quad \tilde{\mathbf{g}} = \begin{bmatrix} \mathbf{g}(t_1) \\ \mathbf{g}(t_2) \\ \vdots \\ \mathbf{g}(t_M) \end{bmatrix}. \tag{25}$$

Thus, we have

$$\tilde{\mathbf{x}}_{n+1} = \tilde{\mathbf{x}}_n + \left( \tilde{\mathbf{J}}_n \tilde{\mathbf{H}} - \tilde{\mathbf{P}} \right) \left( \tilde{\mathbf{Q}} \tilde{\mathbf{x}}_n - \tilde{\mathbf{g}}_n \right), \tag{26a}$$

where $\tilde{\mathbf{H}}$, $\tilde{\mathbf{P}}$, and $\tilde{\mathbf{Q}}$ are constant matrices. Note that subscript $n$ is attached to $\tilde{\mathbf{J}}$ and $\tilde{\mathbf{g}}$ to indicate that they are changing with $\tilde{\mathbf{x}}_n$. In some problems where the magnitude of nonlinear acceleration is not very large, as will be shown in subsection 3.7 using low-Earth-orbit as an example, the Jacobian matrix $\tilde{\mathbf{J}}_n$ can be approximated by a constant matrix. In that case, we can use an even simpler iteration formula instead of Eq. (26a).

$$\tilde{\mathbf{x}}_{n+1} = \tilde{\mathbf{x}}_n + \tilde{\mathbf{J}}_{\mathbf{HP}} \left( \tilde{\mathbf{Q}} \tilde{\mathbf{x}}_n - \tilde{\mathbf{g}}_n \right), \tag{26b}$$

where $\tilde{\mathbf{J}}_{HP} = \tilde{\mathbf{J}}\tilde{\mathbf{H}} - \tilde{\mathbf{P}}$ is a constant matrix.

Finally, Eq. (26a) or Eq. (26b) can be used to iteratively solve for the values of $\mathbf{x}(t)$ at collocated time points. We set the basis functions as orthogonal polynomials. Herein, the first kind of Chebyshev polynomials are adopted, and the collocation points in each time interval are selected as Chebyshev-Gauss-Lobatto (CGL) nodes. The details for $\tilde{\mathbf{H}}$, $\tilde{\mathbf{P}}$, and $\tilde{\mathbf{Q}}$ are provided in appendix.

A flow chart of LVIM method is provided in Fig. 1. A similar flowchart has been introduced in our previous work [29]. However, a more concise and more efficient iteration formula, as is derived in Eq. (26), is applied herein.

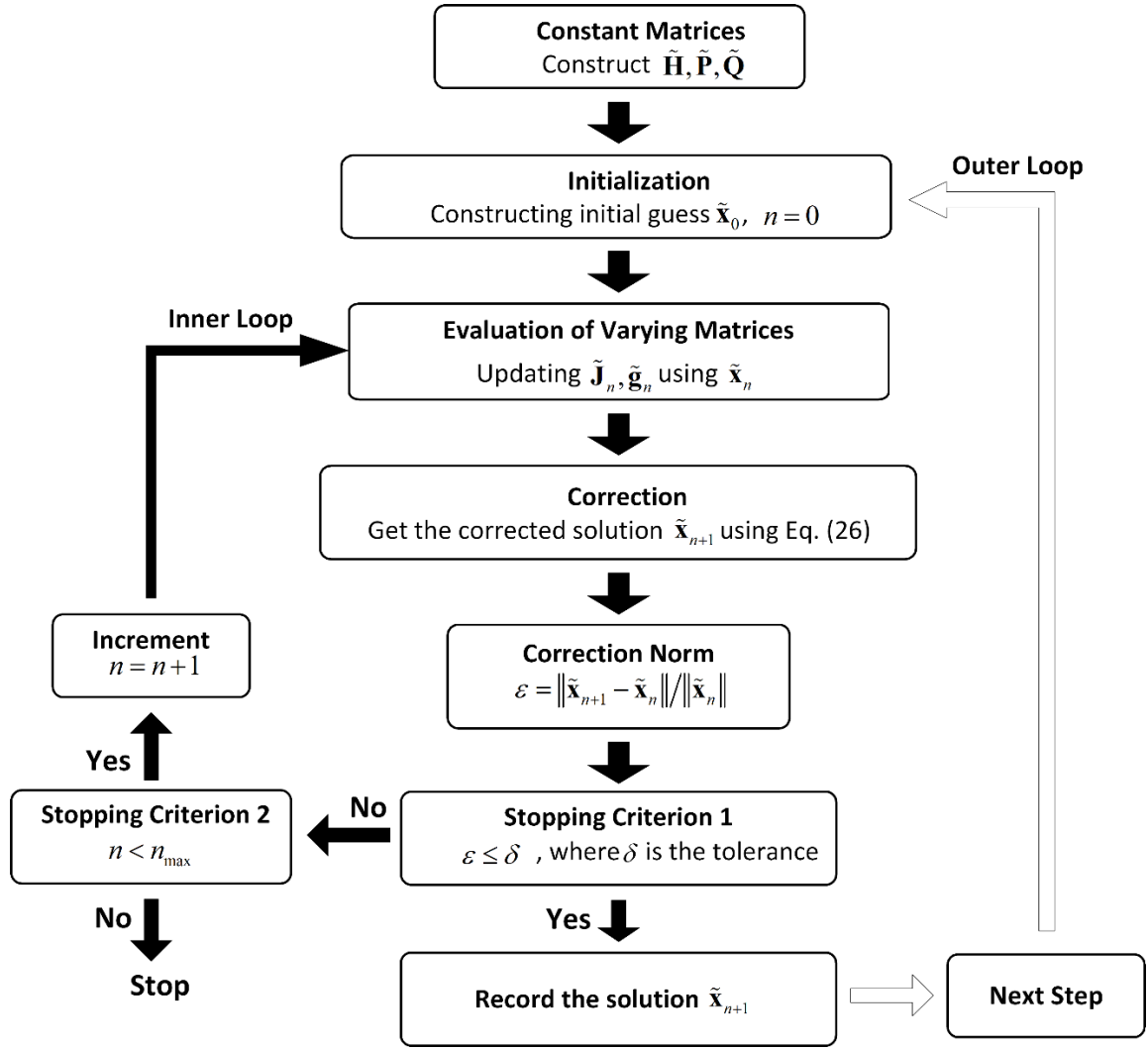

Figure 1. The flow chart of LVIM for solving nonlinear differential equations

## 3. Numerical Results and Discussions

The LVIM is applied herein to several renowned low-dimensional nonlinear problems in classical physics, astrophysics, fluid mechanics, solid mechanics and orbital mechanics. The numerical simulation is carried out in MATLAB R2016b, using a Dell laptop with Intel Core i7-7500U CPU. Note that no GPU acceleration and parallel processing is adopted in the following examples. The computation is conducted only on CPU for consistent comparison, because *ode45* itself have no parallel processing features. However, the performance of LVIM on GPU can be roughly assessed through the total number of iterations needed in the computation.

## 3.1 Blasius Equation

In the problem of laminar boundary layer over a flat plate, the Navier-Stokes equation and the boundary conditions are

$$\frac{\partial u}{\partial x} + \frac{\partial v}{\partial y} = 0$$

$$u\frac{\partial u}{\partial x} + v\frac{\partial u}{\partial y} = \nu \frac{\partial^2 u}{\partial y^2} \tag{27}$$

$$y = 0,\ u = v = 0;\ y = \infty,\ u = U,$$

By introducing the self-similar variables

$$\eta = \frac{y}{\delta(x)} = y\sqrt{\frac{U}{2\nu x}},\ \psi = \sqrt{2\nu U x} f(\eta), \tag{28}$$

and

$$u(x,y) = \frac{\partial \psi}{\partial y} = U f'(\eta),\ v(x,y) = -\frac{\partial \psi}{\partial x} = \frac{1}{2}\sqrt{\frac{\nu U}{x}}\left[\eta f'(\eta) - f(\eta)\right], \tag{29}$$

the Blasius equation is obtained as

$$2f''' + f f'' = 0,\ f = f(\eta), \tag{30}$$

subjected to the boundary conditions

$$\eta = 0,\ f = 0,\ f' = 0;\ \eta = \infty,\ f' = 1 \tag{31}$$

This boundary value problem can be transformed into a pair of initial value problems by transformation groups [30]. The pair of initial value problems are

$$2F''' + F F'' = 0,\ F = F(\xi),$$

$$\xi = 0,\ F = 0,\ F' = 0,\ F'' = 1, \tag{32}$$

and

$$2f''' + f f'' = 0,\ f = f(\eta),$$

$$\eta = 0, \ f = 0, \ f\,' = 0, \ f\,'' = [F\,'(\infty)]^{-3/2}. \tag{33}$$

Eqs. (32, 33) can be solved successively using LVIM. The results are presented in Fig. 2 and compared with those obtained by *ode45* function in MATLAB. The configurations of these two methods (LVIM and *ode45*) and the computational time are listed in Table 1.

Table 1 Configurations and computational time of LVIM and *ode45* in solving Blasius equation

| **Methods** | **Configurations** | **Computational Time** |
|---|---|---|
| LVIM | $N = 5$, $\Delta t = 0.5$, $tol = 10^{-10}$ | 0.008s |
| *ode45* | $RelTol = 10^{-12}$, $AbsTol = 10^{-15}$ | 0.06s |

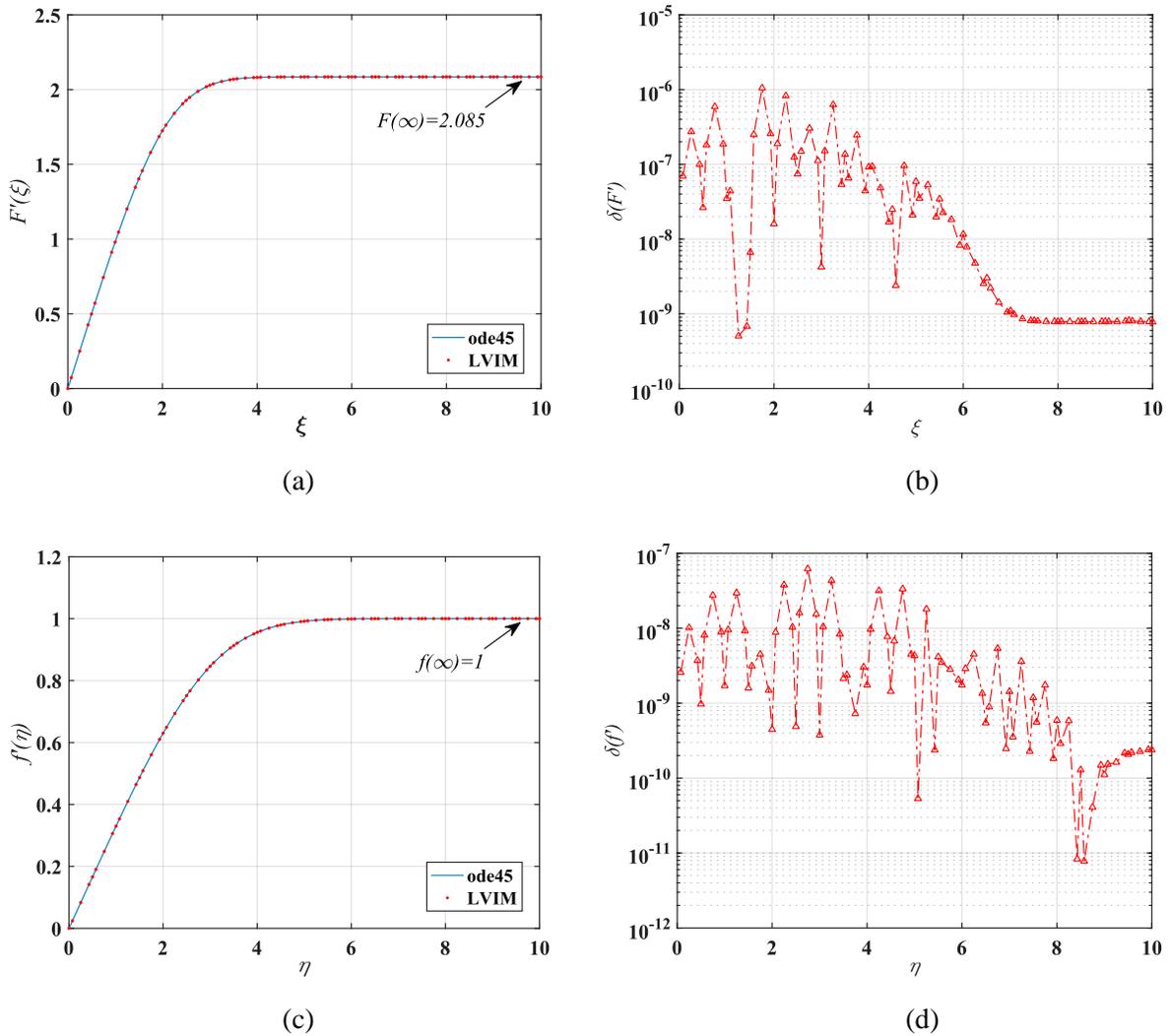

Figure 2. Comparison of the numerical results between LVIM and *ode45* in MATLAB. Left panel: $F'(\xi)$ and $f'(\eta)$ obtained by LVIM and *ode45*. Right panel: computational discrepancy between LVIM and *ode45*.

The normalized stream function obtained by LVIM and *ode45* shows that the flow velocity increases rapidly as $\eta$ increases towards the freestream. The flow velocity is almost the same as the free stream at $\eta = 6$ according to the results. For larger viscosity $\nu$, it can be seen from Eqs. (28, 29) that the flow velocity increases slower in the direction towards free stream, thus the laminar boundary layer becomes thicker. The same thing happens to larger value of $x$, as the flow continues from the leading edge, the thickness of boundary layer also increases.

The results in Fig. 2 shows that the LVIM and *ode45* agrees very well in solving this problem. The right panel of Fig. 2 indicates that the discrepancy between the numerical solutions of these two methods are very small (less than $10^{-6}$). However, it is noted that the computational speed of LVIM is almost 10 times higher than *ode45*, even though *ode45* has been well coded for its optimal performance as a built-in function of MATLAB.

### *3.2 Emden-Chandrasekhar and Chandrasekhar's White Dwarf Equation*

The Emden-Chandrasekhar equation is

$$\frac{1}{\xi^2}\frac{d}{d\xi}\left(\xi^2\frac{d\psi}{d\xi}\right) = e^{-\psi},$$

$$\xi = 0, \ \psi = 0, \ \frac{d\psi}{d\xi} = 0. \tag{34}$$

It describes the density distribution of a spherically symmetric isothermal gas sphere subjected to its own gravitational force, where $\xi$ is the dimensionless radius and $\psi$ is related to the density of the gas sphere.

The Chandrasekhar's white dwarf equation is

$$\frac{1}{\eta^2}\frac{d}{d\eta}\left(\eta^2\frac{d\phi}{d\eta}\right) + (\phi^2 - C)^{3/2} = 0,$$

$$\phi(0) = 1, \ \phi'(0) = 0. \tag{35}$$

where $\phi$ measures the density of white dwarf, and $\eta$ is the dimensionless radius. $C$ is a constant related with the density of white dwarf at the center, which is set as $C = 0.3$ herein. These two problems are solved with LVIM and *ode45*. The results are presented in Figs. 3 and 4.

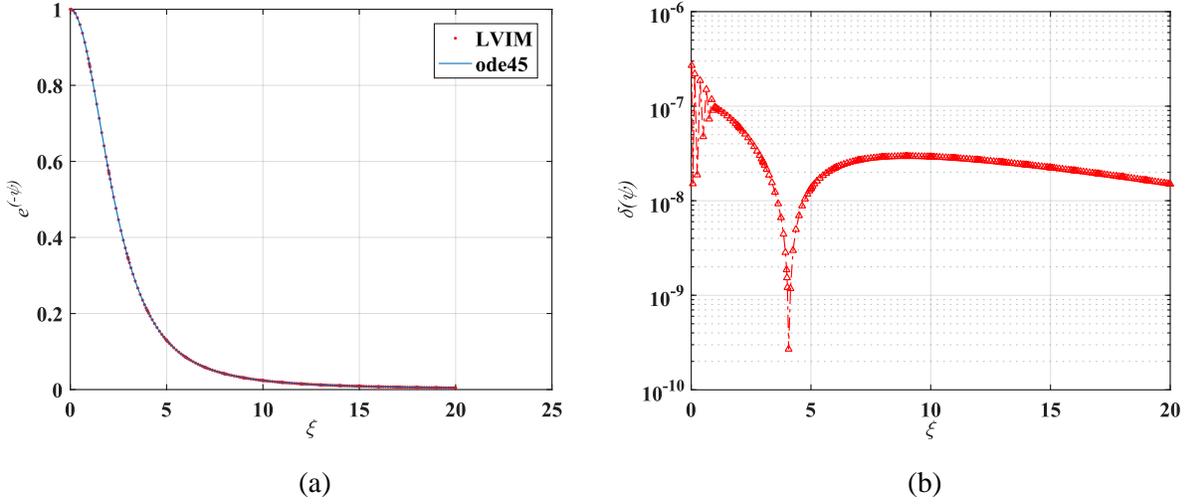

Figure 3. Comparison of the numerical results between LVIM and *ode45* in MATLAB. (a) $e^{-\psi}$ obtained by LVIM and *ode45*, (b) computational discrepancy between LVIM and *ode45*.

Table 2 Configurations and computational time of LVIM and *ode45* in solving Emden-Chandrasekhar equation

| Methods | Configurations | Computational Time |
|---|---|---|
| LVIM | $N=13$, $\Delta t=1$, $tol=10^{-10}$ | 0.006s |
| *ode45* | $RelTol=10^{-12}$, $AbsTol=10^{-15}$ | 0.05s |

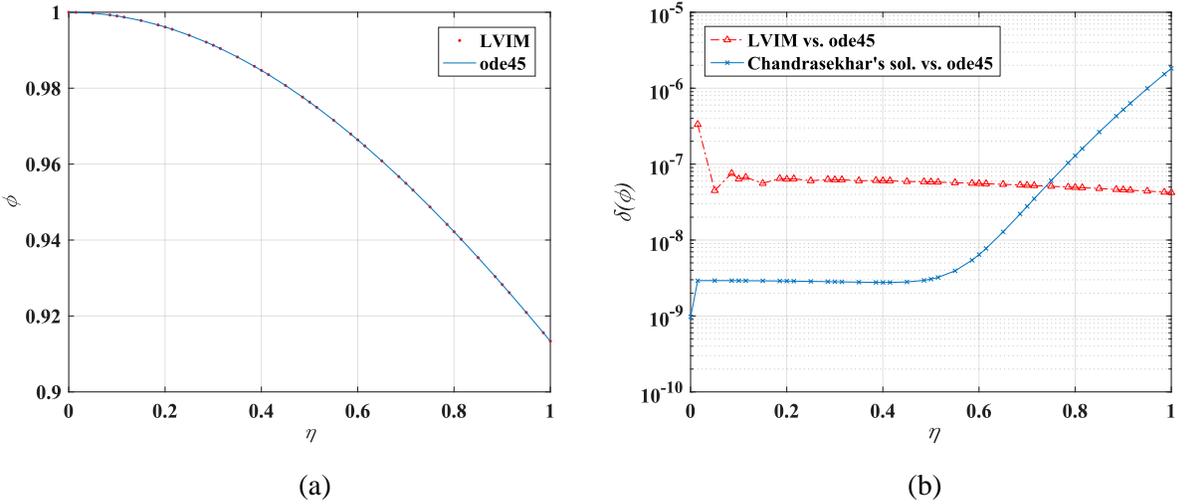

Figure 4. Comparison of the numerical results between LVIM and *ode45* in MATLAB. (a) $\phi$ obtained by LVIM and *ode45*, (b) computational discrepancy between LVIM and *ode45*.

Table 3 Configurations and computational time of LVIM and *ode45* in solving Chandrasekhar's white dwarf equation

| Methods | Configurations | Computational Time |
|---------|----------------|--------------------|
| LVIM | $N=5$, $\Delta t=0.1$, $tol=10^{-10}$ | 0.004s |
| ode45 | $RelTol=10^{-12}$, $AbsTol=10^{-15}$ | 0.03s |

The results show that LVIM can solve these two problems accurately and with very low computational cost. The computational speed of LVIM is almost 10 times higher than *ode45* with the computational discrepancy lower than $10^{-6}$. It is also found that the accuracy of Chandrasekhar's solution deteriorates significantly after $\eta = 0.6$.

### 3.3 Mathieu Equation

The Mathieu equation is

$$\frac{d^2 x}{dt^2} + (\delta - \varepsilon \cos t) x = 0. \tag{36}$$

where the parameters are set as $\delta = 0.5$, $\varepsilon = 0.1$. The initial condition is set as

$$t = 0, \ x = 1, \ \frac{dx}{dt} = 0. \tag{37}$$

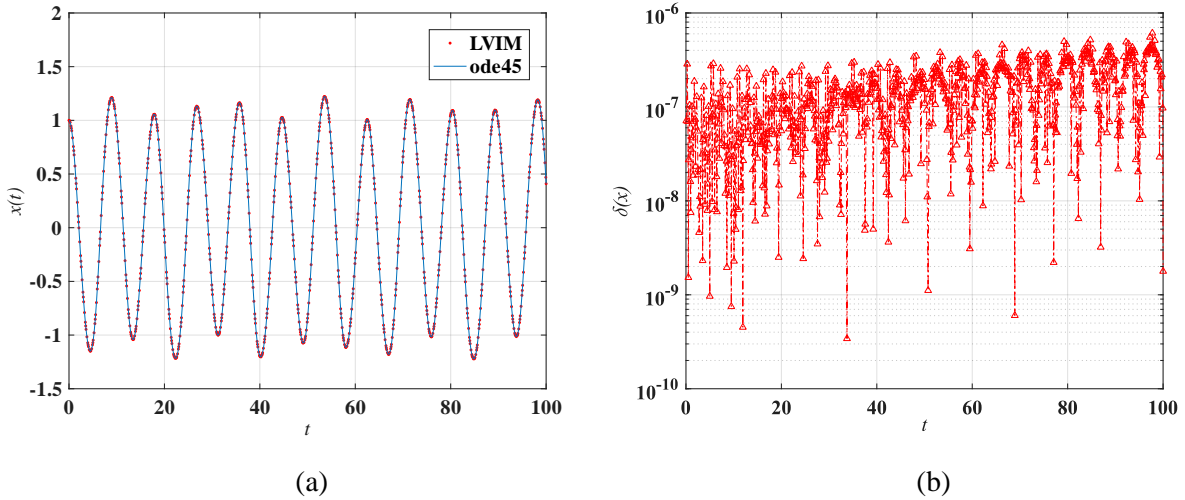

(a)            (b)

Figure 5. Comparison of the numerical results between LVIM and *ode45* in MATLAB. (a) $x(t)$ obtained by LVIM and *ode45*, (b) computational discrepancy between LVIM and *ode45*.

Table 4 Configurations and computational time of LVIM and *ode45* in solving Mathieu equation

| Methods | Configurations | Computational Time |
|---------|----------------|--------------------|
| LVIM | $N=5$, $\Delta t=0.5$, $tol=10^{-10}$ | 0.06s |
| ode45 | $RelTol=10^{-12}$, $AbsTol=10^{-15}$ | 0.24s |

The Mathieu equation is useful in studying nonlinear vibration problems with periodic forcing and the periodic motion of nonlinear autonomous systems. However, since nonlinearity in real physical systems is diminished in Mathieu equation, it may generate some unrealistic solutions such as unbounded motion, because of resonances between the forcing frequency and the oscillator's natural frequency. The stability of the solution of Mathieu equation is dependent on the parameters $\delta$ and $\varepsilon$. The set of parameters used herein generates stable quasi-periodic motion. However, if using parameters such as $\delta = 0.5$, and $\varepsilon = 1$, the resulted motion will become unbounded. It should be noted that in real physics, when the resonance causes the amplitude of the motion to increase, the nonlinear relation between period and amplitude will detune the resonance, thus avoid unbounded motion.

Through numerical simulation using LVIM and *ode45*, a stable quasi-periodic solution is obtained in Fig. 5. The computational discrepancy between these two methods is less than $10^{-6}$. The computational time of LVIM is almost a quarter of that of *ode45*.

### *3.4 Pendulum Equation*

The pendulum equation is

$$\frac{d^2\theta}{dt^2} + \frac{g}{l}\sin\theta = 0 \tag{38}$$

where $g/l$ is set as 1 for convenience. The initial condition is set as

$$t = 0, \ \theta = 3.1329, \ \frac{d\theta}{dt} = 0 \tag{39}$$

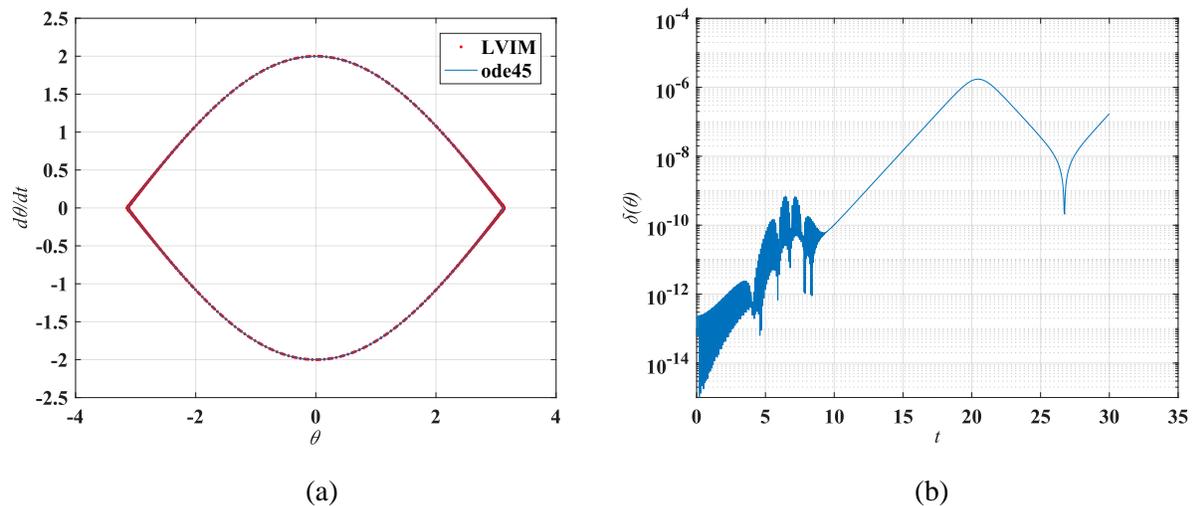

(a)          (b)

Figure 6. Comparison of the numerical results between LVIM and *ode45* in MATLAB. (a) $d\theta/dt$ obtained by LVIM and *ode45*, (b) computational discrepancy $\delta(\theta)$ between LVIM and *ode45*.

Table 5 Configurations and computational time of LVIM and *ode45* in solving pendulum equation

| Methods | Configurations | Computational Time |
|---------|----------------|--------------------|
| LVIM | $N=5$, $\Delta t = 0.1$, $tol = 10^{-10}$ | 0.03s |
| ode45 | $RelTol = 10^{-12}$, $AbsTol = 10^{-15}$ | 0.13s |

In the case of small amplitude, the pendulum motion is nearly sinusoidal, and the frequency is 1. However, for oscillation with large amplitude, the periodic motion is no longer sinusoidal, as shown in Fig. 6 (a). The relationship between amplitude and frequency of a simple pendulum is nonlinear. The frequency changes slowly at first when the amplitude starts increasing from zero, but it changes much faster for relatively larger amplitude. For very large amplitude, the frequency approaches to zero since the period becomes very long. The change of frequency with respect to the amplitude is depicted in Fig. 7.

Using LVIM and *ode45*, the numerical results are obtained in Fig. 6. The configurations of these two methods and the computational time are recorded in Table 5. It is shown that the computational discrepancy is smaller than $10^{-6}$, while the computational time of LVIM is much less than that of ode45.

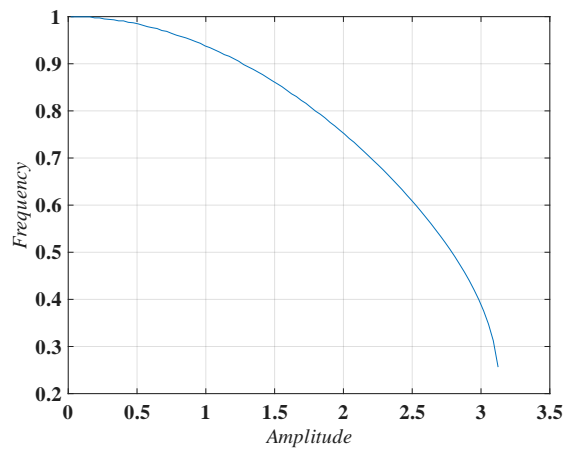

Figure 7. Amplitude-frequency curve of pendulum

### *3.5 Buckled Bar*

The equation of a buckled bar under vertical dead force is written as

$$\text{Vertical dead load: } EI\frac{d^2\theta}{ds^2} = -P\sin\theta, \tag{40}$$

where $E$, $I$, $\theta$, $s$, $P$ are modulus of elasticity, area moment of inertial, bending angle, arc length, and acting force respectively.

Similarly, the equations of a buckled bar under perpendicular follower force and tangent follower force are written as

Perpendicular follower load: $EI\dfrac{d^2\theta}{ds^2} = -P\cos(\theta-\alpha)\sin\theta$,  (41)

Tangent follower load: $EI\dfrac{d^2\theta}{ds^2} = -P\sin(\theta-\alpha)\sin\theta$,  (42)

where $\alpha$ is the bending angle at the top end of the bar. All the three cases are subjected to the boundary condition:

$$s=0,\ \theta=0;\ s=1,\ d\theta/ds=0.\qquad(43)$$

These problems are solved using LVIM in conjunction with shooting method. The results are presented in Fig. 8.

In all three cases, it is found that a slight increase of load above critical can produce considerable deflection of the bar. For large load such as $P=25$, the deflection curves of buckled bar under different type of loading become very different. For relatively small load, the vertical dead loading leads to larger deflection angle, while for large load, the deflection angle caused by the follower force is larger than that by dead force. As the external load becomes very large, saying $P=50$, there could exist multiple solution as shown in Fig. 8 (a).

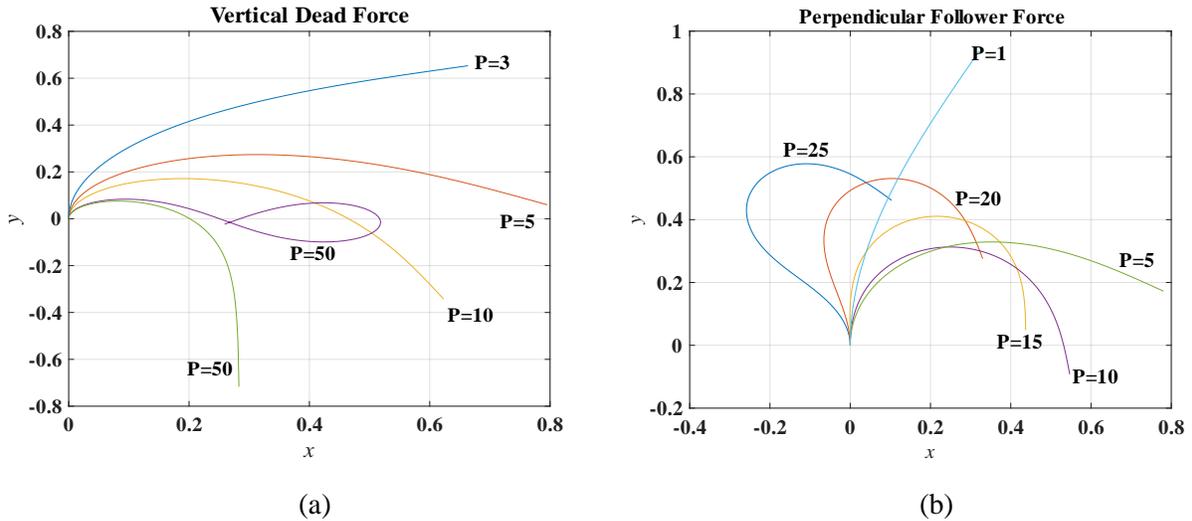

(a)  (b)

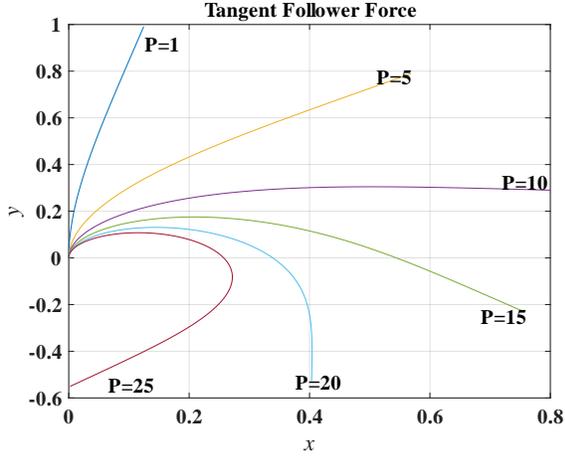
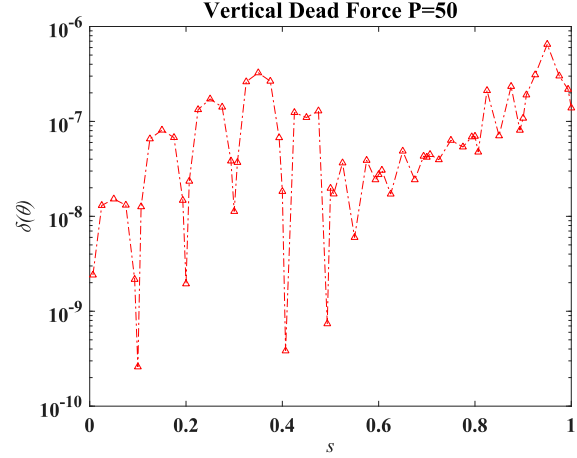

(c)          (d)

Figure 8. Deflection curves of a buckled bar under different types of load. (a) vertical dead loading, (b) perpendicular follower loading, (c) tangent follower loading, (d) computational discrepancy $\delta(\theta)$ between LVIM and *ode45*.

Table 6 Configurations and computational time of LVIM and *ode45* in solving buckled bar equation

| Methods | Configurations | Computational Time |
|---|---|---|
| LVIM | $N = 7$, $\Delta s = 0.1$, $tol = 10^{-10}$ | 0.005s |
| ode45 | $RelTol = 10^{-12}$, $AbsTol = 10^{-15}$ | 0.04s |

Both LVIM and *ode45* are used to solve the buckled bar problem with vertical dead load $P = 50$ through shooting method. The result shows that the computational discrepancy between these two methods is less than $10^{-6}$. In each shooting process, the computational time of LVIM is around 0.005s, while ode45 takes around 0.04s.

### 3.6 Elastica

The Euler's derivation of elastica is

$$\frac{dy}{dx} = \frac{a^2 - c^2 + x^2}{\sqrt{(c^2 - x^2)(2a^2 - c^2 + x^2)}}, \tag{44}$$

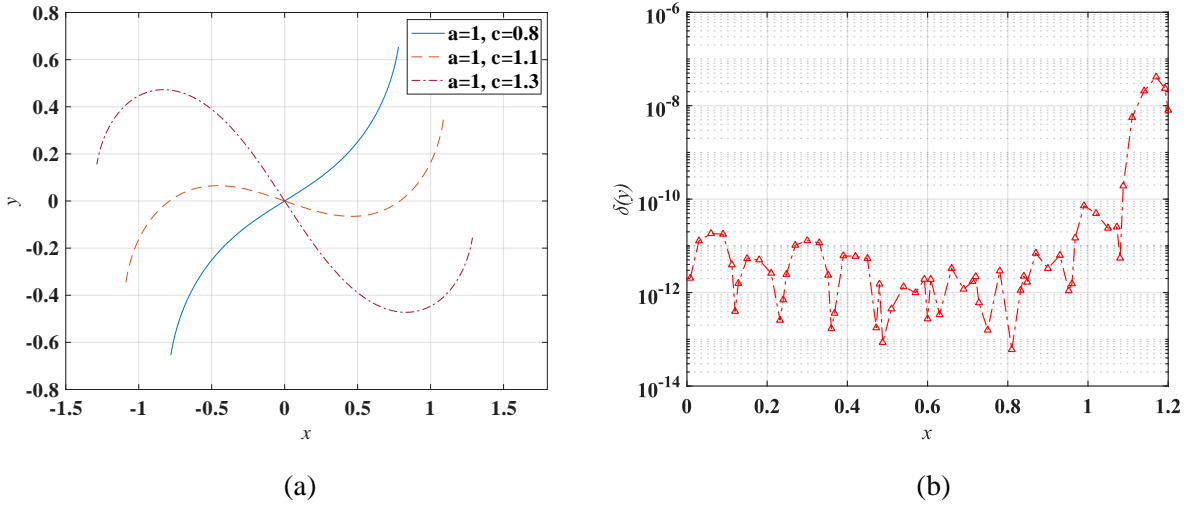

(a)  (b)

Figure 9. Deflection curves of a buckled bar under different types of load. (a) elastica figures for different parameters $a$ and $c$, (b) computational discrepancy $\delta(y)$ between LVIM and *ode45*.

Table 7 Configurations and computational time of LVIM and *ode45* in solving Euler's elastica equation

| Methods | Configurations | Computational Time |
|---|---|---|
| LVIM | $N=7$, $\Delta x = 0.12$, $tol = 10^{-10}$ | 0.003s |
| *ode45* | $RelTol = 10^{-12}$, $AbsTol = 10^{-15}$ | 0.008s |

Although there could be numerous elastic curves, they can be classified into several different cases according to Euler's parameters. The numerical examples we provided herein fall in the categories of $0 < c < a$, $a < c < a\sqrt{1.651868}$, and $a\sqrt{1.651868} < c < a\sqrt{2}$ respectively. Among all the cases, there are three special elasticas for certain parameters including $c=0$, $a=0$, and $c=a$, corresponding to a straight line, a circle, and a rectangular elastica.

Since Euler's elastica equation is one dimension less than the other problems (after transformation into first order differential equations) we solved above, the computational speed of *ode45* improves slightly. However, in this problem the computational time of LVIM is still much less than *ode45*, while the computational error of LVIM is less than $10^{-6}$.

## 3.7 Low-Earth-Orbit Tracking Problem

The governing equation of a low-Earth-orbit (LEO) object is

$$m_i \frac{d^2 \mathbf{q}_i}{dt^2} = -\frac{\partial U}{\partial \mathbf{q}_i}, \tag{45}$$

where $\mathbf{q} = [x, y, z]^T$ is the position vector, $U$ is the potential function of Earth gravity field. It is noted that the 70 deg Earth Gravity Model (EGM) 2008 is used herein, where 4900 nonlinear terms need be calculated in order to accurately evaluate the potential function $U$. The initial condition is set as

$$t = 0, \ x = -0.3889e6, \ y = 7.7388e6, \ z = 0.6736e6,$$

$$dx/dt = -3.5794e3, \ dy/dt = 0, \ dz/dt = 6.1997e3. \tag{46}$$

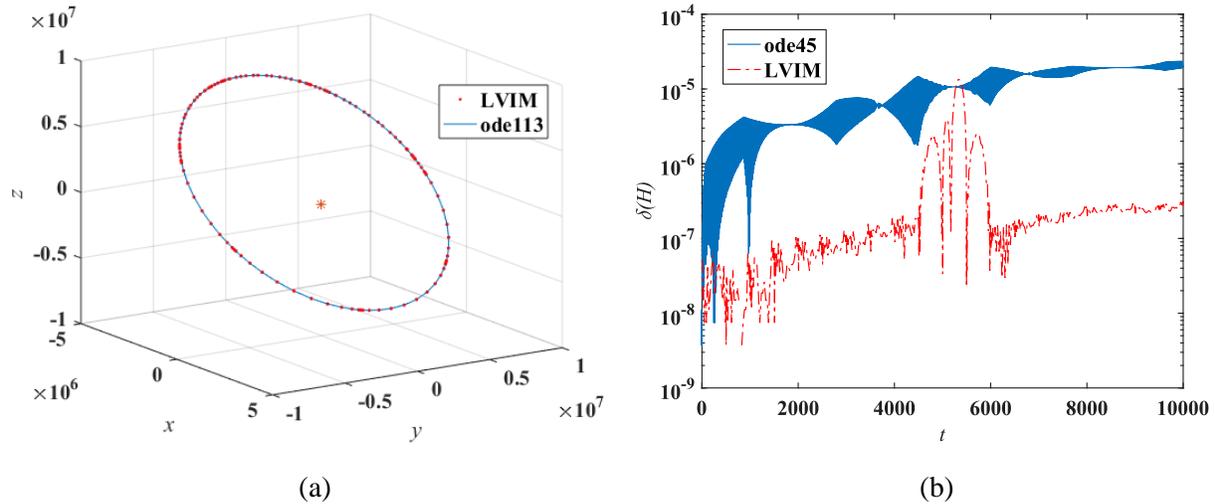

(a)                                          (b)

Figure 10. Comparison of the numerical results between LVIM and *ode45* in MATLAB. (a) trajectory of the target body obtained by LVIM and *ode45*, (b) computational discrepancy $\delta(x)$ between LVIM and *ode45*.

Table 8 Configurations and computational time of LVIM and *ode45* in solving four-body problem

| Methods | Configurations | Num. of Iterations | Computational Time |
|---|---|---|---|
| LVIM | $N = 26$, $\Delta t = 500$, $tol = 10^{-8}$ | 163 | 30s |
| ode45 | $RelTol = 10^{-12}$, $AbsTol = 10^{-15}$ | 13812 | 150s |

The success of space missions highly relies on the computational accuracy of orbital motion. In the case of low-Earth orbit, the long-term trajectory of spacecraft could be affected significantly by a very small perturbation term of gravity force. To achieve a relatively high accuracy, a high order gravitational model should be used. In that case, most computational time is spent in evaluating the gravitational force terms. To save computational time and improve the efficiency, one should reduce the total number of evaluations of force terms in solving the dynamical system.

In Table 8, it is shown that the total number of force evaluations is only 163 by using LVIM, which is almost hundred times less than that using *ode45*. The total computational time of LVIM is 30s, much less than that of *ode45* (150s). Note that for a consistent comparison, the element-wise operations of LVIM are merely conducted on CPU without parallel computation, thus in each round of evaluation, the computational time of LVIM is longer than *ode45*. However, if the element-wise operations of LVIM are parallelly conducted on GPU, the computational time for each round of evaluation could be reduced to the same as that of *ode45*. In that case, LVIM could be almost two magnitudes faster than *ode45*.

In terms of accuracy, the computational error of Hamiltonian of the LEO two body system is recorded in Fig. 10. It shows that the accuracy of LVIM is much higher than *ode45*, although the latter one has been set to its highest threshold of accuracy.

## 4. Conclusion

The proposed LVIM and its related algorithm are used to solve seven representative nonlinear problems in different disciplines of nonlinear science. The first six problems are low-dimensional and have simple nonlinear terms. The numerical results show that LVIM can solve these problems with high accuracy and efficiency. Compared with the built-in *ode45* function of MATLAB, the computational speed of LVIM is almost one magnitude higher, even without GPU acceleration and parallel processing. The last problem we explored is a low-Earth-orbit problem, where the gravity force is evaluated by almost 5000 nonlinear terms. In such a complex problem, the LVIM achieves much higher accuracy than *ode45*, even the latter one has been set to its highest threshold of accuracy. Moreover, the computational time of LVIM is only one fifth of *ode45* when the algorithms are carried out on CPU. According to the total number of iterations needed by LVIM, it is approximated that its computational efficiency can be further improved by twenty times with GPU acceleration. In all, the numerical results verified that the proposed LVIM can be applied extensively in various problems of nonlinear science, and its related algorithm is very accurate and efficient in solving even very complicated nonlinear problems with a large number of degrees of freedom. We propose to include the present LVIM algorithms in an open source code. The application of LVIM in molecular dynamics and high-dimensional nonlinear structural problems will be presented in our future work.

## Appendix

In Eq. (24), $\tilde{\mathbf{H}} = \tilde{\mathbf{P}}\tilde{\mathbf{T}} - \tilde{\mathbf{T}}\tilde{\mathbf{P}}$. Note that $\tilde{\mathbf{P}}$ and $\tilde{\mathbf{Q}}$ are abbreviations for $\tilde{\mathbf{R}}\tilde{\mathbf{P}}\tilde{\mathbf{R}}^{-1}$ and $\tilde{\mathbf{R}}\tilde{\mathbf{Q}}\tilde{\mathbf{R}}^{-1}$ respectively.

$$\tilde{\mathbf{P}} = \begin{bmatrix} \mathbf{P} & & & \\ & \mathbf{P} & & \\ & & \ddots & \\ & & & \mathbf{P} \end{bmatrix}, \tilde{\mathbf{Q}} = \begin{bmatrix} \mathbf{Q} & & & \\ & \mathbf{Q} & & \\ & & \ddots & \\ & & & \mathbf{Q} \end{bmatrix}, \tag{A.1}$$

$$\mathbf{P} = \begin{bmatrix} \int_{t_0}^{t_1} \mathbf{\Phi}(\tau)d\tau \\ \int_{t_0}^{t_2} \mathbf{\Phi}(\tau)d\tau \\ \vdots \\ \int_{t_0}^{t_M} \mathbf{\Phi}(\tau)d\tau \end{bmatrix} \begin{bmatrix} \mathbf{\Phi}(t_1) \\ \mathbf{\Phi}(t_2) \\ \vdots \\ \mathbf{\Phi}(t_M) \end{bmatrix}^{-1}, \mathbf{Q} = \begin{bmatrix} \dot{\mathbf{\Phi}}(t_1) \\ \dot{\mathbf{\Phi}}(t_2) \\ \vdots \\ \dot{\mathbf{\Phi}}(t_M) \end{bmatrix} \begin{bmatrix} \mathbf{\Phi}(t_1) \\ \mathbf{\Phi}(t_2) \\ \vdots \\ \mathbf{\Phi}(t_M) \end{bmatrix}^{-1}, \tilde{\mathbf{T}} = \begin{bmatrix} t_1\mathbf{I} & & & \\ & t_2\mathbf{I} & & \\ & & \ddots & \\ & & & t_M\mathbf{I} \end{bmatrix}. \tag{A.2}$$

$\tilde{\mathbf{R}}$ is a row rearranging matrix that moves the $(dM+m)th$ row of unit matrix to the $(mD+d)th$ row, where $d=1,2,...,D$, $m=1,2,...,M$. $D$ is the dimension of governing system. $M$ is the number of collocation nodes.

In this paper, $\mathbf{\Phi}(t)=[\phi_1(t),\phi_2(t),...,\phi_N(t)]$ is selected as Chebyshev polynomials of the first kind. It is generated by $\phi_{nb}(t)=\cos((nb-1)\arccos t)$, $nb=1,2,...,N$. Thus, we have

Table A.1 Expressions of basis functions (Chebyshev polynomials of the first kind)

| | |
|---|---|
| $\mathbf{\Phi}(t)$ | $\phi_1(t)=1$, $\phi_2(t)=\cos(\arccos t)$, ..., $\phi_N(t)=\cos((N-1)\arccos t)$ |
| $\dot{\mathbf{\Phi}}(t)$ | $\dot{\phi}_1(t)=0$, $\dot{\phi}_2(t)=1$, ..., $\dot{\phi}_N(t)=\sin\left((N-1)^2\arccos(t)/\sqrt{1-t^2}\right)$ <br> (If $t=\pm 1$, let $\dot{\phi}_N(-1)=(-1)^N N^2$, $\dot{\phi}_N(1)=N^2$) |
| $\int_{t_0}^{t}\mathbf{\Phi}(\tau)d\tau$ | $\int_{-1}^{t}\phi_1(\tau)d\tau=t+1$, $\int_{-1}^{t}\phi_2(\tau)d\tau=t^2/2-1/2$, ..., <br> $\int_{-1}^{t}\phi_N(\tau)d\tau=\left\{\dfrac{\cos(N\arccos t)-\cos(N\arccos(-1))}{N}-\dfrac{\cos((N-2)\arccos t)-\cos((N-2)\arccos(-1))}{N-2}\right\}\bigg/2$ |

Noting that we are only using Chebyshev polynomials that are defined in $|t|\le 1$, the real time should be rescaled to fit in this domain.

## Acknowledgements

We thank Texas Tech University for its supports on this research.